\documentclass[secnum,leqno]{rims-bessatsu}

\usepackage[all]{xy}

\usepackage{amssymb,amsmath,color}

\newtheorem{theorem}{Theorem}

\newtheorem{conj}[theorem]{Conjecture}

\newtheorem{cor}[theorem]{Corollary}

\newtheorem{prop}[theorem]{Proposition}
\theoremstyle{remark}

\newtheorem{example}[theorem]{Example}

\def\l{{\mathfrak l}}

\def\O{{\omega}}

\newcommand{\BZ}{{\mathbb Z}}

\newcommand{\G}{{\Gamma}}
\newcommand{\SL}{{SL_2(\BZ)}}

 \let\a\alpha  \let\b\beta  \let\c\gamma  
  
  \let\l\lambda   
      
\let\GL\Lambda
\def\C{\mathbb C}
\def\G{\mathbf G}

\def\sgn{\rm sgn}

\def\GL{\mathbf{GL}}
\def \GL2 {{\text{GL}_2}}

\def\Gal{{\rm Gal}}
\def\Frob{{\rm Fr}}
\def\F{{\mathbb F}}
  
\def\Z{{\mathbb Z}}

\def\Q{{\mathbb Q}}

\def\O{\mathcal O}

\def\G{\Gamma}

\def\ord{{\mathrm ord}}

\catcode`,\active

\catcode`\,12

\title{Atkin and Swinnerton-Dyer congruences and noncongruence modular forms }
\TitleHead{Atkin and Swinnerton-Dyer congruences}
\author{Wen-Ching Winnie \textsc{Li}\footnote{Department of Mathematics, Pennsylvania State University,
University Park, PA 16802, USA ~and National Center for Theoretical Sciences, Mathematics Division,
National Tsing Hua University, Hsinchu 30013, Taiwan, R.O.C.\newline e-mail: \texttt{wli@math.psu.edu}} ~and Ling \textsc{Long}
\footnote{Department of Mathematics, Louisiana State University, Baton Rouge, LA 70803 USA \newline e-mail: \texttt{llong@lsu.edu}}}
\AuthorHead{Wen-Ching Winnie \textsc{Li} and Ling \textsc{Long}}
\classification{11F11, 11F33, 33C20}

\support{The first author is supported by the NSF grant DMS-1101368.  The second author is supported by NSF grant DMS-1001332. Part of the material was based on the talk given by the first author at RIMS on Dec. 5, 2012. This paper was written while the second author was visiting Cornell University as an AWM Michler scholar. She is in debt to both AWM and Cornell University. The computations were made using the mathematical software {\sf Sage}.}
\begin{document}
\maketitle
\begin{abstract}
  Atkin and Swinnerton-Dyer congruences are special congruence recursions satisfied by coefficients of noncongruence modular forms. These are in some sense $p$-adic analogues of Hecke recursion satisfied by classic Hecke eigenforms. They actually appeared in different contexts and sometimes can be obtained using the theory of formal groups. In this survey paper, we introduce the Atkin and Swinnerton-Dyer congruences, and discuss some recent progress on this topic.
\end{abstract}

\section{Introduction}
Atkin and Swinnerton-Dyer (ASD) congruences in the title refer to the congruences of the form
\begin{eqnarray}\label{ASDgeneral}
 a_{np} - A_p a_n + \mu_p p^{k-1}a_{n/p} \equiv 0 \mod{p^{(k-1)(1 + \ord_p n)}}  \quad \rm{for~all}~ n \ge 1, \end{eqnarray}where $p$ is a prime, $a_n$ are integral over $\Z_p$, $A_p$ is an algebraic integer,  and $\mu_p$ is a root of unity. The congruence means that the left hand side divided by the modulus is integral over $\Z_p$.
 As usual, $a_x = 0$ if $x$ is not an integer. The original purpose was to study the arithmetic properties of the Fourier coefficients $a_n= a_n(f)$ of a weight-$k$ cusp form $f$ for a finite index subgroup $\G$ of $\SL$. When $\G$ is a congruence subgroup and $f$ is an eigenfunction of the Hecke operator at $p$, then with $A_p$ being the eigenvalue and congruence replaced by equality, this is the familiar three term recursive {Hecke} relation. When $\G$ is a noncongruence subgroup, which is the majority, such $f$ may be regarded as playing the role of a Hecke eigenfunction at $p$. This remarkable observation was made by Atkin and Swinnerton-Dyer in their seminal paper \cite{A-SD}, which initiated a systematic study of the arithmetic of noncongruence modular forms. Assume that the modular curve of $\G$ has a model defined over $\Q$ such that the cusp at infinity is $\Q$-rational, $k \ge 2$ and the space $S_k(\G)$ of weight-$k$ cusp forms for $\G$ is $d$-dimensional. Based on their numerical evidence, Atkin and Swinnerton-Dyer expected  $S_k(\G)$ to contain a basis with this congruence property for good primes $p$,  although the basis would vary with $p$.  Using the {$2d$-dimensional $\ell$-adic Galois} representations attached to $S_k(\G)$,
Scholl in \cite{Sch85b} proved that such a basis exists at primes $p$ where the action of the Frobenius has $d$ distinct $p$-adic unit eigenvalues. In  \S6.3, we give some examples of $S_k(\G)$ for which the ASD expectation holds {for almost all primes},
and also an example where the 3-term ASD congruence holds for only half of the primes. In the paper \cite{Sch85b} Scholl also showed that all forms $f$ in $S_k(\G)$ with $p$-adically integral Fourier coefficients satisfy a similar, but longer, $(2d+1)$-term congruence for all $n\ge 1$:
\begin{eqnarray*}
a_{np^d}(f)+A_1a_{np^{d-1}}(f)+\cdots+A_d a_n(f)+\cdots+ A_{2d}a_{n/p^d}(f)\equiv 0 \mod p^{(k-1)(1+\ord_pn)}.
\end{eqnarray*} Here $T^{2d} + A_1T^{2d-1} + \cdots + A_{2d} \in \Z[T]$ is the characteristic polynomial of the geometric Frobenius at $p$ under the Galois representation.  See \S2 for details.

Recently, Kazalicki and Scholl \cite{KS13} extended the above congruence to include weakly holomorphic exact weight-$k$ {cusp} forms for $\G$, but with weaker modulus $p^{(k-1)\ord_pn}$ instead. This is discussed in \S \ref{ss:7}. In \S \ref{ss:app}, we describe an application of the ASD congruences to an open conjecture characterizing genuine noncongruence modular forms.

The congruences discovered by Atkin and Swinnerton-Dyer actually also appeared in different contexts. For instance, when the modulus is $p^{1+\ord_pn}$ (with $k=2$), the ASD-type congruences can be obtained using the theory of formal groups, recalled in \S \ref{ss:2}. These congruences are closely related to Dwork's work on $p$-adic hypergeometric series, see \S \ref{ss:3}. Section \ref{ss:4} is devoted to some geometric backgrounds for ASD congruences, including $p$-adic analogues of the Selberg-Chowla formula and solutions of certain ordinary differential equations. We also exhibit some known and some conjectural supercongruences, which are ASD-type congruences that are stronger than what can be predicted from the formal group laws.  We end the paper by discussing another type of congruences discovered by Atkin, which {led to} the development of $p$-adic modular forms.
\smallskip

{\noindent \textsc{Acknowledgment.} The authors would like to thank the referee whose
comments greatly improve the exposition of this article.}

\section{Congruence and noncongruence modular forms}\label{ss:mf}
\subsection{Finite index subgroups of $\SL$ and modular forms}
{It is well-known that all finite index subgroups of $SL_n(\Z)$ for $n \ge 3$ are congruence subgroups (cf. \cite{BLS, BMS}).} This property no longer holds for
 $\SL$. In fact, among its finite index subgroups, noncongruence ones out number the congruence ones.
One way to see this is from the defining fields of the associated modular curves. More precisely, the modular curve $X_\G$ of a finite index subgroup $\G$ of $\SL$ is a Riemann surface
obtained as the orbit space of the action of $\G$ on the Poincar\'e upper half-plane  {$\frak H$} via fractional linear transformations, compactified by adding finitely many cusps. It is known to have a model defined over a number field. When $\G$ is a congruence subgroup $\G_0(N)$ or $\G_1(N)$, the modular curve is defined over $\Q$, and for the principal congruence subgroup $\G(N)$, its modular curve has a model defined over the cyclotomic field $\Q(e^{2\pi i /N})$. On the other hand,
a celebrated result of Bely{\u\i} asserts that

\begin{theorem}[Bely{\u\i}, \cite{Bel}] A smooth irreducible projective curve defined over a number field is isomorphic to a modular curve $X_\G$ for (infinitely many) finite index subgroup(s) $\G$ of $\SL$.
\end{theorem} Thus simply by considering curves defined over number fields one sees easily that $\SL$ contains far more noncongruence subgroups than congruence ones.

Given a finite index subgroup $\G$ of $\SL$, recall that a weight-$k$ modular form for $\G$ is a holomorphic function $f$ on the upper half-plane satisfying
\begin{equation}\label{modular}
 f(z) = (cz+d)^{-k} f \left (\frac{az+b}{cz+d}\right ) \quad \text{for~all}~ \begin{pmatrix}
  a&b\\c&d \end{pmatrix}\in \G {\text{~and~all}~ z\in \frak H}
  \end{equation} and the extra condition that $f$ is holomorphic at the cusps of $\G$. It is called a \emph{cusp form} if it vanishes at all cusps of $\G$. We call it a \emph{congruence form} if it is for a congruence subgroup; otherwise it is called a \emph{noncongruence form}. Denote by $M_k(\G)$ (resp. $S_k(\G)$) the space of all weight-$k$ modular forms (resp. cusp forms) for $\G$.

When $k=2$, the forms $f \in S_2(\G)$ may be identified with the holomorphic differential 1-forms {$f(z)dz$} on $X_\G$. The de Rham space $H^1(X_\G,\C)$ is $2g$-dimensional, where $g$ is the genus of $X_\G$, {spanned by} holomorphic and anti-holomorphic 1-forms on $X_\G$; these two spaces are dual to each other with respect to the cup product on $H^1(X_\G,\C)$.

\subsection{Congruence modular forms} The arithmetic for congruence forms is well-understood, after  being studied for over one century. The two main ingredients are the Hecke theory and {$\ell$-adic} Galois representations. The newform theory says that it suffices to study the arithmetic of newforms of weight $k$, level $N$, and character $\chi$ \cite{AL,Miyake, Li}. Such a form $f$ is a common eigenfunction of the Hecke operators at the primes $p$ not dividing $N$. The Fourier coefficients $a_n(f)$ of $f$ with the leading coefficient $a_1(f) = 1$ are algebraic integers satisfying
$$ a_{mn}(f) = a_m(f)a_n(f) \qquad \text{for}~(m, n) = 1, $$ and the $3$-term recursive relation
\begin{equation}\label{hecke}
 a_{np}(f) - a_p(f)a_n(f) + \chi(p)p^{k-1}a_{n/p}(f) = 0 \qquad \text{for~all}~ p\nmid N~\text{and~all}~ n \ge 1.
\end{equation}The work of Eichler-Shimura \cite{shim1} (for $k = 2$) and Deligne \cite{Del} (for $k \ge 3)$ associates to $f$ a compatible family of 2-dimensional {$\ell$-adic} representations $\rho_{\ell,f}$ of the absolute Galois group $G_\Q$ over $\Q$, unramified outside $\ell N$ such that the characteristic polynomial of $\rho_{\ell, f}$ at the Frobenius at $p \nmid \ell N$ is $H_p(T) = T^2 - a_p(f)T + \chi(p)p^{k-1}$. {Since $\rho_{\ell,f}$ is of motivic nature, by the Weil conjecture proved by Deligne, we have} $|a_p(f)| \le 2 p^{(k-1)/2}$ for all $p\nmid N$, {which is the celebrated Ramanujan-Petersson conjecture.} The above developments originated from Ramanujan's observations on the discriminant Delta function {in $S_{12}(\SL)$}  $$\Delta(z):=\eta(z)^{24}=q\prod_{n\ge 1}(1-q^n)^{24}, \quad \text{where}~q=e^{2\pi i z}.$$

\subsection{Noncongruence modular forms}
In comparison,  modular forms for noncongruence subgroups are far more mysterious than their congruence counterparts due to the lack of effective Hecke operators, which was conjectured by Atkin and proved by Serre \cite{Thompson80} for noncongruence subgroups normal in $\SL$ and Berger \cite{berger94} in general. More precisely, if one mimics what's done for congruence forms by defining a Hecke operator at $p$ for a noncongruence subgroup $\G$ by using the $\G$-double coset represented by a $2 \times 2$ matrix with entries in $\Z$ and determinant $p$, then, as shown in \cite{berger94}, this operator is the composition of the trace map from $\G$ to $\G^c$, the smallest congruence subgroup containing $\G$, followed by the usual Hecke operator at $p$ on $\G^c$. Unfortunately the trace map annihilates all genuine noncongruence forms, hence no information can be drawn for the noncongruence forms we are interested in. However, it is not difficult to construct  noncongruence modular forms as long as we keep an eye on the ramification. For weight 0 meromorphic modular forms, namely \emph{modular functions}, we have
 \begin{theorem}[Atkin and Swinnerton-Dyer, \cite{A-SD}]\label{thm:ASD1}
   An algebraic function of the (modular) $j$-function is a modular function if and only if, as a function of  $j$-function,  it only ramifies at 3 points  1728, 0 and infinity with ramification indices 2, 3, and arbitrary, respectively.
 \end{theorem}

A systematic investigation on noncongruence modular forms was initiated by the work of Atkin and Swinnerton-Dyer \cite{A-SD}. They showed some similarities between congruence and noncongruence forms.
Given a finite index subgroup $\G$ of $\SL$, for convenience, assume that the modular curve $X_\G$ has a model defined over $\Q$ under which the cusp $\infty$ is a $\Q$-rational point. Then there is an integer $M$, {divisible by the widths of the cusps of $\G$ and the primes $p$ where the   cusps are no longer distinct under reduction of $X_\G$ at $p$, }
such that $S_k(\G)$ contains a basis whose  Fourier coefficients are integral over $\Z[\frac{1}{M}]$, see \cite{Sch85b}. Suppose $k \ge 2$ is even and $S_k(\G)$ has dimension $d$. From their numerical examples, Atkin and Swinnerton-Dyer observed that for good primes{\footnote {In their paper \cite{A-SD} Atkin and Swinnerton-Dyer called these primes good without giving any definition. It turns out that in all numerical examples where the ASD congruences hold, the {$A_p(i)$} is a $p$-adic unit. Scholl has shown that if $\rho_\ell(\Frob_p)$ is semi-simple with half of the eigenvalues being $p$-adic units, then $p$ is a good prime. } } $p\nmid M$ and $k$ even,  the space $S_k(\G)$ possesses a basis $\{f_i \}_{1 \le i \le d}$ with $p$-adically integral Fourier coefficients $a_n(f_i)$ and for each $i$ there exists an algebraic integer $A_p(i)$ with $|A_p(i)|\le 2p^{(k-1)/2}$ such that
\begin{eqnarray}\label{ASDcongruence}
 a_{np}(f_i)-A_p(i)a_n(f_i)+p^{k-1}a_{n/p}(f_i)\equiv 0 \mod p^{(k-1)(1+\ord_pn)}, \quad \forall n\ge 1.
\end{eqnarray}
The striking similarity between this and \eqref{hecke} suggests that the $p$-adic theory for noncongruence modular forms could be fruitful as well.

To understand what Atkin and Swinnerton-Dyer discovered,  Scholl in \cite{Sch85b} constructed {a compatible family of {$2d$-dimensional} $\ell$-adic representations $\rho_\ell = \rho_{\G, k, \ell}$ of the absolute Galois group $G_\Q$} attached to the space $S_k(\G)$ for $k > 2$.   For each prime $\ell$, $\rho_\ell$ is an extension of the construction of Deligne, obtained from the first \'etale cohomology of $X_\G(\overline \Q)$ with coefficients in the sheaf which is the $(k-2)$nd  symmetric power of some local system. Scholl proved that $\rho_\ell$ is  unramified outside $\ell M$ such that for any prime $p \nmid \ell M$
the characteristic polynomial  $H_p(T)= T^{2d} + A_1T^{2d-1} + \cdots + A_{2d}$ of the geometric Frobenius
under $\rho_\ell$ is in $\Z[T]$  with all roots of the same absolute value $p^{(k-1)/2}$. In fact, $H_p(T)$ can be computed by explicit formulas of counting points on elliptic curves  over finite fields. Moreover, by a comparison theorem between \'etale and crystalline cohomologies, Scholl proved that for $p>k-2$ and $p\nmid M$, any $f \in S_k(\G)$ with $p$-adically integral Fourier coefficients $a_n(f)$ satisfies the congruence that for all $n\ge 1$
\begin{multline}\label{schollcongruence}
a_{np^d}(f)+\cdots+A_d a_n(f)+\cdots+ A_{2d}a_{n/p^d}(f)  \equiv 0  \mod p^{(k-1)(1+\ord_pn)}.
\end{multline}  In particular,
when $H_p(T)$ has $d$ distinct $p$-adic unit roots, the above long congruence can be reduced to 3-term congruences satisfied by {the Fourier coefficients of} a basis of $S_k(\G)$, as observed by Atkin and Swinnerton-Dyer.

Scholl representations are of motivic nature. According to the Langlands philosophy, the $L$-functions of Scholl representations should coincide with the $L$-functions of automorphic representations of certain reductive groups.
The congruences \eqref{ASDcongruence} and \eqref{schollcongruence} can then be interpreted as congruence relations between Fourier coefficients of noncongruence forms and those of automorphic forms. For 2-dimensional Scholl representations of $G_\Q$, their modularity follows from a renowned conjecture of Serre established by Khare and Wintenberger \cite{Khare-Wintenberger09} and various modularity lifting theorems  \cite[et al.]{SW2, SW, DM03, dml, kisin9, ALLL}.
What makes Scholl representations interesting is that they cannot be decomposed into {$2$-dimensional} pieces in general and hence provide a fertile testing ground for Langlands philosophy. For more details in this regard, see a survey \cite{Li12}. In \S \ref{ss:app}, we will give an application of the ASD congruences as well as the automorphy of Scholl's representations \cite{LL12}. In the proof, we also used the best known bound for the coefficients $a_n(f)$ of noncongruence cusp forms obtained by Selberg, {namely $|a_n(f)| = O(n^{k/2-1/5})$.}  {In \cite{selberg65} Selberg gave an example to show that the Ramanujan-Petersson conjecture fails for noncongruence cusp forms}.

To illustrate the above results we exhibit an example below.  Let  $\G^1(5)$  be the group consisting of matrices in $\SL$ which become lower triangular unipotent when modulo 5. It has 4 cusps,  $\infty$, 0, $-2$, and $-5/2$, and admits a normalizer $A=\begin{pmatrix}
  -2&-5\\1&2
\end{pmatrix}$ which swaps cusps $\infty$ and {$-2$}. Let $E_1$, $E_2$ be weight-3 Eisenstein series of $\G^1(5)$, which have simple zeros at all cusps except $\infty$ and $-2$, respectively, and nonvanishing elsewhere. Both of them have integer Fourier coefficients:
 \begin{eqnarray*}
  E_1(z)& =& 1- 2q^{1/5} - 6q^{2/5} + 7q^{3/5} + 26q^{4/5} + \cdots,\\
    E_2(z)& =& q^{1/5}- 7q^{2/5} + 19q^{3/5} -23q^{4/5} + \cdots.
\end{eqnarray*}For more details, see \cite[\S 4]{LLY05}. Thus $t=\frac{E_2}{E_1}$ is a modular function for $\G^1(5)$ with a simple zero and a simple pole, both located at the cusps. By Theorem \ref{thm:ASD1}, $\sqrt t$ is a modular function for an index-2 subgroup $\G_2$ of $\G^1(5)$. Also, $$f= E_1\sqrt{t}=\sqrt{E_1E_2}=\sum_{n\ge 1} a_n(f)q^{n/10}=q^{1/10}-\frac{3^2}2q^{3/10}+\frac{3^3}{2^3}q^{5/10}+\frac{3\cdot 7^2}{2^4}q^{7/10}+\cdots$$ is a weight-3 cusp form for $\G_2$, which in fact generates $S_3(\G_2)$. The group $\G_2$  is a noncongruence subgroup because the Fourier coefficients of $\sqrt{E_1E_2}$ have unbounded denominators. No congruence cusp form behaves this way, since it is a linear combination  of Hecke eigenforms, whose Fourier coefficients are algebraic integers. This distinction is conjectured to be a criterion to distinguish congruence forms from genuine noncongruence forms with algebraic Fourier coefficients. This conjecture is of fundamental importance and is very useful in many applications.  It will be discussed later in \S7. Nevertheless, {if the coefficients of a noncongruence modular form lie in a number field, then they are integral at all but finitely many places.}

For $S_3(\G_2)$  the  ASD congruences proved by Scholl assert that, for all primes $p>3$, there are $A_p,B_p$ in $\Z$ such that
$$a_{np^r}(f)-A_pa_{np^{r-1}}(f)+B_pa_{np^{r-2}}(f)\equiv 0 \mod p^{2r}, \quad \forall n,r\ge 1.$$ In particular, the corresponding 2-dimensional  {$\ell$-adic} Scholl representation of $G_\Q$ is reducible when restricted to $\Gal(\overline \Q/\Q(\sqrt{-1}))$ due to the finite order linear operator  {on the representation space induced from the matrix $A$ above. On the space of forms, this operator sends $E_1$ to $E_2$ and $E_2$ to $-E_1$.}  Consequently, one can check that $A_p$ agrees with the $p$th coefficient of the weight-3 Hecke eigenform $\eta(4z)^6=q\prod_{n\ge 1} (1-q^{4n})^6$ and $B_p=\left ( \frac{-1}p\right )p^2$, where $\left ( \frac{-1}p\right )$ is the Legendre symbol. For more details, see \cite{LLY05}.

\section{ASD congruences and 1-dimensional commutative formal group laws}\label{ss:2}
\subsection{ASD congruences for elliptic curves}
The inspiration of the congruence \eqref{ASDcongruence} observed by Atkin and Swinnerton-Dyer  came from what they proved for weight-$2$ cusp forms for $\G$ {such that the genus of $X_\G$ is one}. In this case the modular curve $X_\G$ is an elliptic curve $E: \ y^2=x^3+Ax+B$ with $A, B \in \Z$. Let $p>3$ be a prime such that $E$ has good reduction modulo $p$.  Let $\xi$ be a local uniformizer  of $E$ {at the point at infinity} 
which is either $-\frac xy$ or $-\frac xy$ plus higher order terms with coefficients in $\Z$. The coefficients of the holomorphic differential $1$-form $\frac{dx}{2y}=\sum_{n\ge1} a_n \xi^n\frac{d\xi}{\xi}$ satisfy
\begin{equation}\label{eq:ASD-EC}
a_{np^r}-(p+1-\#[E/\F_p])a_{np^{r-1}}+pa_{np^{r-2}}\equiv 0 \mod p^r, \quad \forall n,r\ge 1.
\end{equation} This can be explained by using 1-dimensional commutative formal group law, which we will recall in the next section by following the development in \cite{Kibelbek}. For more information on $d$-dimensional commutative formal group laws, the reader is referred to  \cite{Haz, Ditters, Kibelbek}.

\subsection{1-dimensional commutative formal group law (1-CFGL)}
 A 1-CFGL over a characteristic $0$ commutative integral domain $R$ with units, such as $\Z$, $\Z[1/M]$, or $\Z_p$, is a formal power series $G(x,y)$ in  $R[[x,y]]$  satisfying
\begin{itemize}
\item  $G(x,y)=x+y+\sum_{i,j\ge 1} c_{i,j}x^iy^j, c_{i,j}\in R$,
\item (Associativity) $G(x,G(y,z))=G(G(x,y),z)$,
\item (Commutativity) $G(x,y)=G(y,x)$.
\end{itemize}

Corresponding to each formal group is a unique normalized \emph{invariant differential}
{$$\omega(y)=[\partial_x G(0,y)]^{-1}dy=(1+\sum_{i\ge 2} a_i y^{i-1}) dy, a_i\in R$$ where $\partial_x$ means partial derivative with respect to the first variable.
The corresponding \emph{strict formal logarithm} is defined to be $$\ell(y):=\int \omega= y+\frac{a_2}{2}y^2+\frac{a_3}{3}y^3+\cdots$$ so that $G(x,y)=\ell^{-1}(\ell(x)+\ell(y))$.}

Two 1-CFGLs $G(x,y)$ and $\bar G(x,y)$ over $R$ are said to be \emph{isomorphic} if there exists a formal power series $\phi(x)=ax+\text{higher terms}$ in $R[[x]]$ with $a\in R^\times$  such that $\phi(G(x,y))=\bar G(\phi(x),\phi(y))$. If $a=1$, then the isomorphism is said to be \emph{strict}. Thus, over the field of fractions of $R$, each 1-CFGL is strictly isomorphic to the additive CFGL $G(x,y)=x+y$ via the formal logarithm.
 \begin{prop}\label{prop:changevariable-CFGL}
  If $f(x)=\sum_{n\ge 1} \frac{a_n}{n}x^n$ is the strict formal logarithm of a 1-CFGL $G(x,y)$ over $R$, then for any $\phi(x)=x+\text{higher terms} \in R[[x]]$, $f(\phi(x))=\sum_{n\ge 1} \frac{b_n}nx^n$ is the strict formal logarithm of a 1-CFGL, which is strictly isomorphic to $G(x,y)$.
\end{prop}

{Denote by $S$ the set of formal power series of the form $f=\sum_{n\ge 1} \frac{a_n}n x^n$ with $a_n\in R$. It is an $R$-module. On $S$ we define the following operators with integers $m \ge 1$ and $\lambda\in R$:
\begin{itemize}
  \item (\emph{Frobenius}) $F_mf(x)=\sum_{n\ge 1}\frac{a_{mn}}nx^n$,
  \item (\emph{Verschiebung}) $V_mf(x)=\sum_{n\ge 1}\frac{a_{n}}nx^{mn}$,
  \item (\emph{Witt operator}) $[\l]f(x)=\sum_{n\ge 1} \frac{a_n}n\l^nx^n$.
\end{itemize}Note that $[\l]f(x)$ corresponds to the change of variable $x \mapsto \l x$. Thus any formal group isomorphism can be decomposed into a strict isomorphism followed by a Witt operator.

The following theorem characterizes the strict formal logarithms of  1-CFGLs.

\begin{theorem}
$f \in S$ is the strict formal logarithm of a 1-CFGL over $R$ if and only if for  each prime $p$ there exist $\l_{p, i} \in R$ such that
$$F_pf=\sum_{i\ge 1}V_i[\l_{p, i}]f.$$
 \end{theorem}

For the remaining discussion, assume that for any maximal ideal $\wp$ of $R$ with  residual  characteristic $p$, the completion of the localization of $R$ at $\wp$ is an unramified ring extension  of $\Z_p$. Let  $\sigma_\wp$ be the ring automorphism of $R$ sending $r \in R$ to $r^{\sigma_\wp} \equiv r^p \mod \wp$.
On the submodule $S_p$ of $S$ of formal power series of the form $g= \sum_{i\ge 0}\frac{a_{p^i}}{p^i}x^{p^i}$, for each $\mu\in R$ we define the \emph{Hilbert operator} $\{\mu\}$: it sends $g$ to
$\{\mu\}g=\sum_{i\ge 0} \frac{a_{p^i}}{p^i}\mu^{\sigma_\wp^i}x^{p^i}.$
Then the Witt operators on $S_p$ are generated by $F_p, V_p$ and the Hilbert operators.

\begin{theorem}
Suppose $f =\sum_{n\ge 1} \frac{a_n}{n} x^n  \in S$ is the strict formal logarithm of a 1-CFGL over $R$ as above.  Then  there are unique $\mu_{p, i} \in R$ such that
  $$F_pf_{(p)}=\sum_{i\ge0}V_{p^i}\{\mu_{p,i}\} f_{(p)}, \quad \text{where}~~f_{(p)}=\sum_{i\ge 0}\frac{a_{p^i}}{p^i}x^{p^i}.$$
	\end{theorem} }

The theorem above asserts the following relation on coefficients of $f_{(p)}$:
$$a_{p^{n+1}}=\sum_{i=0}^np^ia_{p^{n-i}}\mu_{p,i}^{\sigma_\wp^{n-i}}=\sum_{i=0}^n p^i \mu_{p,i} a_{p^{n-i}}^{\sigma_\wp^{i+1}},$$ where the last equality is due to the commutativity relation $\{\mu\}V_p=V_p\{\mu^{\sigma_\wp}\}$. This leads to the following general congruences:
\begin{equation}\label{eq:ASD-CFGL-2}
  a_{mp^{n+1}}\equiv \sum_{i=0}^n p^i \mu_{p,i} a_{mp^{n-i}}^{\sigma_\wp^{i+1}} \mod \wp^{n+1}, \quad \forall~n, m\ge 1.
\end{equation}

When $a_p\not \equiv 0\mod \wp$, i.e. $p$ is ordinary for the 1-CFGL,  then \eqref{eq:ASD-CFGL-2} can be reduced to a 2-term congruence
\begin{equation}\label{eq:2-term}
  a_{mp^{n+1}} \equiv \a_p a_{mp^{n}}^{\sigma_\wp} \mod \wp^{n+1}, \quad \forall~n,m\ge 1,
\end{equation} for a some $\wp$-adic unit $\a_p$.
\begin{theorem}\label{thm:iso-Hilbert}
{Let $f$, $\wp$ and $\mu_{p,i}$ be as in Theorem 5. Let $g(x)=\sum_{n\ge 1}\frac{b_n}{n}x^n$ be the strict formal logarithm of a 1-CFGL over $R$ which is strictly isomorphic to the 1-CFGL of which $f$ is the strict formal logarithm. Then
\begin{equation}\label{eq:ASD-CFGL}
  b_{mp^{n+1}}\equiv \sum_{i=0}^n p^i \mu_{p,i} b_{mp^{n-i}}^{\sigma_\wp^{i+1}} \mod \wp^{n+1}, \quad \forall~n, m\ge 1.
\end{equation}}
\end{theorem}

\subsection{ASD congruences for elliptic curves continued}
To resume our discussion in \S3.1, one can construct a {1-CFGL}  associated to the  elliptic curve $E: y^2=x^3+Ax+B$ as follows. In a neighborhood $V$ of the point at infinity, each point $P=(x,y)$ on $E$ is marked by  $\xi(P)=-\frac{x}{y}$. Given $P_1,P_2$ in $V$, one can compute $P_3=P_1+P_2$ under the group law on $E$ and express $\xi(P_3)$ as a formal power series  in  $\xi(P_1)$ and $\xi(P_2)$ with coefficients in ${  \Z}$. This formal power series satisfies all axioms of the 1-CFGL with $\frac{dx}{2y}$ expanded in powers of $\xi= -\frac xy$ as its invariant differential. {Namely, $\frac{dx}{2y}$ gives rise to a 1-CFGL over ${ \Z}$ which is naturally identified with the infinitesimal group law of the elliptic curve around the point at infinity.
Honda \cite{Honda-FGL-zeta} proved that this formal group is strictly isomorphic to the formal group constructed from the $L$-function $L(E, s) = \sum_{n \ge 1} b_n n^{-s}$ of $E/\Q$. Thus both group laws yield congruence relations with the same $\mu_{p,0}, \mu_{p,1},...$ by Theorem \ref{thm:iso-Hilbert}. It follows from the definition of  $L(E, s)$ as an Euler product that its coefficients $b_n$ satisfy the relation $b_p = p+1-\#[E/\F_p]$ and
$$b_{np} - b_pb_n + pb_{n/p} = 0, \quad \forall~ n \ge 1$$ for all primes $p$ not dividing the conductor $N$ of $E$. For such $p$ we have $\mu_{p,0} = b_p$, $\mu_{p,1} = -1$ and $\mu_{p,j} = 0$ for $j \ge 2$ so that \eqref{eq:ASD-EC} holds. }

\section{ASD-type congruences and Dwork's result}\label{ss:3}
\subsection{ASD congruences for Legendre family of elliptic curves}\label{ss:3.1}
For {an integer $r \ge 1$} and $\alpha_i,\beta_i \in \mathbb{C}$ with $\beta_i\notin \Z_{\le 0}$, the hypergeometric series {$_{r}F_{r-1}$} is defined by
\[
{}_rF_{r-1}\left[{\a_1, \ldots, \a_r \atop \b_1, \ldots, \b_{r-1}}; x \right] =
\sum_{k= 0}^{\infty} \frac{(\alpha_1)_k(\alpha_2)_k \cdots (\alpha_{r})_k}{(\beta_1)_k \cdots (\beta_{r-1})_k} \cdot \frac{x^k}{k!},
\]
which converges for $|x|<1$. Here $(a)_k := a\cdot(a+1)\cdots(a+k-1)$ is  the Pochhammer symbol. This series satisfies an order $r$ ordinary differential equation in $x$.  Denote by
${}_rF_{r-1}\left[{\a_1, \ldots, \a_r \atop \b_1, \ldots, \b_{r-1}}; x \right]_n$
the truncation of the series after the $x^n$ term.

Let $E_\l: y^2=x(x-1)(x-\lambda)$ be the Legendre family of elliptic curves. Its  Picard-Fuchs equation  is an {order $2$} hypergeometric differential equation with the unique holomorphic solution near zero (up to scalar multiples) equal to
${}_2F_{1}\left[{\frac 12, \frac 12 \atop 1}; \lambda \right].$
Expand $\omega=\frac{dx}{2\sqrt{x(x-1)(x-\l)}}=\sum_{n\ge 1} \frac{a_n}n\xi^n\frac{d\xi}\xi$ with $\xi=-\frac{x}{y}$. {Then $a_n=0$ for $n$ even}. Using a formula due to Beukers (see \cite[pp. 272]{Ditters}), {one gets that  $a_n$ is the coefficient of $x^{n-1}$ in $(x(x-1)(x-\l))^{\frac{n-1}2}$ for odd $n$}. Thus   \begin{equation} \label{a_k} a_{2k+1} =
{}_2F_{1}\left[{-k, -k \atop 1}; \lambda \right](-1)^k.
\end{equation}As $-k$ is a negative integer, the above hypergeometric series terminates at $x^k$. When $\l\in\Q$, the ASD congruence \eqref{eq:ASD-EC}  for ordinary prime $p$ is equivalent to
\begin{equation}\label{2F1cong}
\frac{{}_2F_{1}\left[{\tfrac{1-p^s}{2}, \tfrac{1-  p^s}{2} \atop 1}; \lambda \right] }{
{{}_2F_{1}\left[{\tfrac{1-p^{s-1}}{2}, \tfrac{1-  p^{s-1}}{2} \atop 1}; \lambda \right]
}}\equiv \left(\frac{-1}{p}\right) \beta_{p,\l}  \mod p^s, \, { \forall s\ge 1}
\end{equation} where $\beta_{p,\l}$ is the unit root of $T^2 -(p+1-\#[E_\l/\F_p])T+p$ and hence only depends on the residue of $\l$ in $\F_p$. For details, see \cite{KLMSY}.

We proceed to compare the above with some results of Dwork, which is a special case of Theorem 2 \cite{Dwork-p-cycle}.

\begin{theorem}[Dwork, \cite{Dwork-p-cycle}]\label{thm:Dwork}Let $p>2$ be a fixed prime and let  $K$ be a complete $p$-adic field with $R$ its ring of integers and $p$ a uniformizer.  Let $ B(n)$ be an $R$-{valued sequence}.
  Let $F(X):=\sum_{n\ge 0}B(n)X^n$. Suppose that for all integers $n\ge 0, m\ge0,$  $s\ge 1,$
  \begin{itemize}
  \item[(1)] $B(0)$ is a unit in $R$;
  \item[(2)] $\frac{B(n)}{B(\left [ \frac np\right])} \in R$;
    \item[(3)] $\frac{B(n+mp^{s+1})}{B(\left [ \frac np\right]+mp^s)}\equiv \frac{B(n)}{B(\left [ \frac np\right])} \mod p^{s+1}$.
\end{itemize} Then
  $$F(X)\sum_{j=mp^s}^{(m+1)p^s-1}B(j)X^{pj}\equiv F(X^p) \sum_{j=mp^{s+1}}^{(m+1)p^{s+1}-1}B(j)X^{j}\mod B(m)p^{s+1}R[[X]].$$
\end{theorem}Intuitively, the  condition (3) is about certain $p$-adic continuity of $\frac{B(n)}{B(\left [ \frac np\right])}$, as a function of $n$. It is satisfied by   the $p$-adic (Morita) Gamma function $\G_p(x)$, which is a {continuous function from $\Z_p$ to  $\Z_p^\times$ such that $\G_p(0)=1$ and} for any $x\in \Z_p$, $\G_p(x+1)/\G_p(x)=-x$ if $p\nmid x$ and $\G_p(x+1)/\G_p(x)=-1$ otherwise. It is known that $\G_p(n+mp^{s+1})\equiv \G_p(n) \mod p^{s+1}$, for more details, see  \cite{Cohen}.

From definition, $$n!=(-1)^{n+1}\G_p(1+n) \cdot p^{[n/p]}[n/p]!.$$
Thus
\begin{equation}\label{eq:binom-2n-n}
  \frac{\binom{2n}{n} }{\binom{2[n/p]}{[n/p]}}=-\frac{\G_p(1+2n)}{\G_p(1+n)^2}s_{n,p},
\end{equation} where $s_{n,p}=1$ if $n-p\cdot [\frac np]< p/2$; and  $s_{n,p}=p(2[\frac np]+1)$ otherwise.

By properties of $\G_p(x)$ and $s_{n,p}$, for any positive integer $k$,  $B(n)=\left (\frac{(\frac 12)_n}{n!} \right )^k$ satisfy all the conditions of  Theorem \ref{thm:Dwork}.
 When $k=2$, Theorem \ref{thm:Dwork} implies that for $\lambda\in \Z_p$ such that $E_\l$ has ordinary reduction modulo $p$, there exists a unit $\a_{p,\l}\in \Z_p$
such that
\begin{equation}\label{eq:Dwork-cong}
\frac{{}_2F_1\left[{\frac12, \frac12 \atop 1}; \l \right]_{p^s-1}}
{{}_2F_1\left[{\frac12, \frac12 \atop 1}; \l^p \right]_{p^{s-1}-1}}
 \equiv \a_{p,\l}  \mod p^s, \quad \text{for all } s\ge 1.
\end{equation}
Unlike $\beta_{p,\l}$ above, $\a_{p,\l}$ may vary if we replace $\l$ by $\l^p$. However, $\left ( \frac{-1}p \right ) \b_{p,\l}=\a_{p,\hat \l}$ where $\hat \l$ is the Teichm\"uller lift of $\l$ modulo $p\Z_p$ to $\Z_p$. In fact, the discrepancy between \eqref{eq:Dwork-cong} and  \eqref{2F1cong} basically lies in the differences between the Witt operators and Hilbert operators introduced earlier. Using the result of Dwork and Ditters \cite{Ditters}, Kibelbek showed that if we let
$A(n) = {}_rF_{r-1}\left[{\frac12, \ldots, \frac12 \atop 1, \ldots, 1}; \lambda \right]_{n-1}$,
 then $\sum \frac{A(n)}nx^n$ is the strict formal logarithm of a 1-CFGL over $\Z[1/M][\l]$ for some integer $M$ \cite{Kibelbek-private}. In \cite{KLMSY},  a geometric interpretation of these 1-CFGLs was given explicitly based on \cite{Stienstra}.  By relaxing the condition  of Theorem \ref{thm:Dwork}, Dwork's result implies that truncated hypergeometric series with rational upper parameters (and lower parameters to be all 1's) give rise to CFGLs, integral at almost all primes, which are not necessarily 1-dimensional. It is natural to ask whether one can find isomorphic formal groups arising from an explicit algebraic equation. In fact, many of them are realized using hypersurfaces in weighted projective spaces, including many of geometric objects being Calabi-Yau manifolds. Meanwhile, for the untruncated hypergeometric series, they correspond to objects like periods, which we will illustrate using the Legendre family below.

\section{ ASD congruences, periods, differential equations, and related topics}\label{ss:4}
\subsection{Periods, Picard-Fuchs equations, and modular forms}\label{ss:4.1}
A period of an elliptic curve $E$ is an integral $\displaystyle \int_{\gamma} \frac{dx}{2y}$ over $\gamma\in H_1(E,\Z)$.
 In general, these are transcendental numbers.  For the Legendre family, the variation of the periods, $p(t)=\int_{\gamma_t} \frac{dx}{2\sqrt{x(x-1)(x-t)}}$, is captured by its Picard-Fuchs (PF) equation alluded to in the previous section. Near 0, the unique (up to scalar multiple) holomorphic solution of this PF equation is
 ${}_2F_{1}\left[{\tfrac12, \tfrac12 \atop 1}; t \right].$
{ There is a choice of the cycles $\gamma_t$ to make $p(t)$ a holomorphic function in $t$
so that $p(t) = C \cdot {}_2F_{1}\left[{\tfrac12, \tfrac12 \atop 1}; t \right]$ for $C$ an algebraic multiple of $\pi$ \cite{BB}.} It is well-known that the Legendre family represents elliptic curves with 2-torsion points, whose moduli space $X(2)$ is parametrized by the classical modular lambda function $\lambda(z):=16\cdot \frac{\eta(2z)^4\eta(z/2)^2}{\eta(z)^6}$.  Setting $t=\l(z)$, one has
\begin{equation}\label{eq:3}
{}_2F_{1}\left[{\tfrac12, \tfrac12 \atop 1}; \lambda(z) \right] = \theta_3^2(z),
\end{equation}  where $\theta_3(z)=\sum_{n\in \Z} q^{n^2/2}$ with $q=e^{2\pi i z}$ is a Jacobi theta function of weight $\frac12$.

\subsection{Complex multiplication and results of Chowla-Selberg}
When the elliptic curve $E$ over a number field has complex multiplication, i.e. its endomorphism ring {$\frak R$} over $\C$ is an order of an imaginary quadratic field $K=\Q(\sqrt{-d})$ with fundamental discriminant $-d$,  all periods are algebraic multiples of a transcendental number $b_K$, depending on $K$. The Selberg-Chowla formula \cite{Chowla-Selberg} describes a choice of $b_K$:
\begin{equation}\label{eq:CS}
  b_K:=\G(\frac12)\prod_{0<a<d}\Gamma(\frac ad)^{\frac{n\varepsilon(a)}{4h}},
\end{equation}
where $n$ is the order of unit group in $K$,
$\varepsilon$ is a primitive quadratic  Dirichlet character  modulo $d$,
and $h$ is the class number of $\Q(\sqrt{-d})$.

\begin{example}\label{eg:period-Q(i)}
$b_{\Q(\sqrt{-4})}=\G(\frac 12)\frac{\G(\frac14)}{\G(\frac 34)}\sim{ \frac{\G(\frac14)^2 }{\G(\frac12)}}$, by the reflection formula $\G(x)\G(1-x)=\frac{\pi}{\sin (\pi x)}.$ { Here $\sim$  means equality} up to an algebraic multiple.\end{example}

Regard the invariant differential $\omega = \frac{dx}{2y}$ of $E$ as an element of $H_{DR}^1(E,\C)$, the dual of $H_1(E,\Z)\otimes _\Z \C$.
The endomorphism ring  $\frak R$ of $E$ over $\C$ acts on the $2$-dimensional space $H_{DR}^1(E,\C)$ with
$\omega$ as a common eigenvector of  $\frak R$.
There is another common eigenvector $\nu$ for  $\frak R$ in  $H_{DR}^1(E, \C)$, defined over $\overline \Q$ and  independent of $\omega$.  Chowla and Selberg further showed that the particular \emph{quasi-period}
 $\int_\gamma \nu$ is an algebraic multiple of  $2\pi \sqrt{-1} / b_K $, hence its transcendental part is also a Gamma quotient. Putting  together, one obtains the relation
\begin{equation}\label{eq:CS2}\int _\gamma \omega \cdot \int_\gamma \nu\sim \pi.
\end{equation}

{From \eqref{eq:3},} this number $b_K$ also plays a role at {singular} values of modular forms. More precisely, for any weakly holomorphic {(i.e. allowing poles at cusps)} modular form $F$ with integral weight $k$ and algebraic Fourier coefficients, it is known that {$F(\tau)\sim (b_K/\pi)^{k}$} for all $\tau\in K$ with $\text{Im}(\tau) > 0$, see {\cite{Zagier-123}}.

\subsection{ASD congruences and Gross-Koblitz formula}
It is worth mentioning that Selberg-Chowla's results have $p$-adic analogues. Their formula \eqref{eq:CS} for $b_K$ inspired the Gross-Koblitz formula \cite{GK}. Let $\pi_p \in \C_p$ be a fixed root of $x^{p-1}+p=0$. Let $\varphi$ be the Teichm\"uller character. The Gross-Koblitz formula states that under a suitable normalization {depending on the choice of $\pi_p$ and the choice of an additive character in the Gauss sum}, we have the Gauss sum
\begin{equation}\label{eq:G-K}g(\varphi^{-j})=-\pi_p^j\G_p(\frac{j}{p-1}), \quad 0\le j\le p-2.
\end{equation} Young gave  another proof of the Gross-Koblitz formula using formal groups constructed from the Fermat curves \cite{Young94}.

\subsection{ASD congruences and $p$-adic periods}\label{ss:5.4} In some sense, ASD congruences also describe  $p$-adic periods.
To illustrate the idea, we use the following example of CM elliptic curve $\displaystyle  E: y^2=x^3+x$, with endomorphism ring  $\frak R=\Z(\sqrt{-1})$ due to the order 4 automorphism $\displaystyle (x,y)\mapsto (-x,\sqrt{-1}y)$. If we expand $\frac{dx}{2y}=\sum_{n \ge 1} a_n\xi^n\frac{d\xi}\xi$ with $\xi=\frac{-x}y$,  we  have  $a_n=\binom{\frac{n-1}2}{\frac{n-1}4}$ when $n\equiv 1 \mod 4$, and $a_n=0$ otherwise by the formula of Beukers \cite[pp. 272]{Ditters} alluded to before.\footnote{Another way to derive a formula for $a_n$ is to use \eqref{a_k} for $E_{-1}: y^2=x(x-1)(x+1)$ which is isomorphic to $E$ over $\Q(\sqrt{-1})$. This approach requires an additional evaluation of Kummer to reach our conclusion.}   For a prime $p\equiv 1 \mod 4$, $E$ is ordinary, and the  ASD congruence is reduced to a 2-term relation:  using \eqref{eq:binom-2n-n} we have, {for all $n\equiv 1 \mod 4$,}
 $$
{\binom{\frac{np^r-1}2}{\frac{np^r-1}4}}\slash{\binom{\frac{np^{r-1}-1}2}{\frac{np^{r-1}-1}4}}=-\frac{\G_p(\frac {1+np^r}2)}{\G_p(\frac{3+np^r}4)^2}= \frac{\G_p(\frac{1-np^r}4)^2}{\G_p(\frac {1-np^r}2)},$$ {where the last equality follows from the $p$-adic analogue of the reflection formula for the Gamma function $\G_p(x)\G_p(1-x)=(-1)^{x_0}$ with $x_0\in\{1,2,\cdots,p\}$ being the residue of $x$ mod $p$.} As $r \to \infty$, we find the limit  $\a_p$ in \eqref{eq:2-term} is {$\a_p=\frac{{\G}_p(
\frac14)^2}{\G_p(\frac12)}$. }This is a $p$-adic analogue of the period computed via Selberg-Chowla formula in Example \ref{eg:period-Q(i)}.

\subsection{An example of weakly holomorphic differentials}Consider the Fermat curve $E: x^3+y^3=1$ with $(x,y)=(0, 1)$ as the origin. Its endomorphism ring is  $\frak R=\Z[\frac{1-\sqrt{-3}}2]$.  Two linearly independent eigenvectors of  $\frak R$ in $H_{DR}^1(E,\C)$ are $\omega=\frac{dx}{y^2}=\sum_{n\ge 1} a_nx^n\frac{dx}x$, a holomorphic differential, and $\nu=\frac{xdx}{y}=\sum_{n\ge 2} b_nx^n\frac{dx}x$, a differential of second kind. At $p\equiv 1 \mod 3$, which are ordinary primes for $E$, we have
$$a_{p^n}=\binom{-\frac23}{\frac{p^n-1}3} \quad \text{and} \quad b_{2p^n}=\binom{-\frac13}{\frac{2(p^n-1)}3}.$$   For details,  see Example 6.3 of \cite{Kat-crys}. By the discussion in \S \ref{ss:2},  $\displaystyle a_{p^{n}}/a_{p^{n-1}}\equiv \a_p \mod p^n$ for all $n\ge 1$, where $\a_p$ can be written in terms of $p$-adic Gamma product.  However, numerical data seem to indicate deeper congruence relations:
\begin{equation}\label{eq:18}
 a_{p^{n}}\equiv \a_p a_{p^{n-1}} \mod p^{2n},\quad b_{2p^{n}} \equiv \frac{p}{\a_p}b_{2p^{n-1}} \mod p^{2n-1}.
\end{equation}

\subsection{Another $p$-adic analogue of Selberg-Chowla relation, Ramanujan-type congruences}
 For the Legendre family, $H^1(E_t,\C)$ is generated by $\omega_t$, a holomorphic differential, and $\partial_t \omega_t$, see \cite{Kat}. Thus, the quasi-period of $E_t$ can be written as a linear combination of ${}_2F_1\left[{\frac12, \frac12 \atop 1}; t \right]$
 and {its derivative with respect to $t$}.  Similar to \eqref{eq:3}, one can relate the quasi-period to  an explicit  modular form of weight $-1$. For details, see \S 3.4 of \cite{CDLNS}.  Using the well-known formula of Clausen:
\begin{equation}\label{eq:clausen}
{}_2F_1\left[{1-c, c \atop 1};x \right]^2 = {}_3F_2\left[{\frac12, 1-c, c \atop 1, 1}; -4x(x-1) \right],
\end{equation}which holds for $c,x\in \C$ such that both hypergeometric series  converge, when $E_t$ has complex multiplication, one can express the formula \eqref{eq:CS2} as
\begin{equation}\label{eq:pi-other}
   \sum_{k=0}^{\infty} \left (\frac{(\frac12)_k}{k!} \right)^3(-4t(t-1))^k(ak+1)=\frac \delta \pi,
\end{equation} for  some computable algebraic numbers  $a,\delta$ depending on $t$. The derivation is given in \cite{CDLNS}. Formulas of this type include
 $$\sum_{k= 0}^{\infty} \left (\frac{(\frac12)_k}{k!} \right)^3(6k+1)\frac{1}{4^k}=\frac 4 \pi$$ by Ramanujan. These so-called Ramanujan-type formulas for $1/\pi$
 were first given by  Borwein-Borwein \cite{BB} and  Chudnovsky-Chudnovsky \cite{CC}. Later, van Hamme discovered several surprising $p$-adic analogues of Ramanujan-type formulas for $1/\pi$ \cite{vanHamme97}.  For instance, he conjectured that for each prime  $p>3$
\begin{eqnarray}
    \sum_{k=0}^{ \frac{p-1}2 } \left (\frac{(\frac12)_k}{k!} \right)^3(6k+1)\frac{ 1}{4^k}&\equiv& \left ( \frac{-1}p \right ) p \mod {p^4}, \label{eq:1.3}
 \end{eqnarray}
   where $\left ( \frac{\cdot }p \right )$ is the Legendre symbol. This conjecture was proved in \cite{Lon11}.
More generally, we have {
\begin{theorem}[Theorem 1 of \cite{CDLNS} with $d=2$]\label{thm:main2}Given $\l\in \overline \Q^{\times}$ such that  $\Q(\l)$ is totally real and $|\l|<1$ for all embeddings,  let $t=\frac{1-\sqrt{1-\l}}2$. Assume that the elliptic curve $E_{t}: y^2=x(x-1)(x-t)$ has complex multiplication. Let $a$ be the algebraic number in \eqref{eq:pi-other} for $E_t$. For each odd prime $p>3$ unramified in $K=\Q(t)$ such that $a,\l$ can be embedded in $\Z_p^\times$ and $E_t$ has good reduction at any prime ideal of the ring of integers of $K$ dividing $(p)$,  we have
 \begin{equation*}
  \sum_{k=0}^{p-1}  \left (\frac{(\frac12)_k}{k!} \right )^3(ak+1)\l^k\equiv  {\sgn} \cdot \left ( \frac{1-\l}p \right ) \cdot p \mod {p^2},
\end{equation*}
where $\left ( \frac{1-\l}p \right )$ is the Legendre symbol of the residue class modulo $p$ of $1-\l$, and $\sgn=\pm 1$,    equal to  $1$ if and only if $p$ is ordinary for $E_{t}$.
\end{theorem}} Meanwhile, numerical data suggest the following ASD-type congruences {for all integers $n\ge1$}, extending the conjecture by Zudilin which is the case $n=1$:
 \begin{equation}\label{eq:22}
		  \sum_{k=0}^{p^n-1} \left (\frac{(\frac12)_k}{k!} \right )^3 (ak+1)\l^k\equiv  {\sgn} \cdot \left ( \frac{1-\l}p  \right )  \cdot p \cdot    \sum_{k=0}^{p^{n-1}-1}  \left (\frac{(\frac12)_k}{k!} \right )^3(ak+1)\l^k  \mod {p^{3n}}.
\end{equation}

\subsection{Supercongruences and complex multiplication}
In \cite{CDE86}, Chowla, Dwork, and Evans gave an improvement of the ASD congruence {related to \S \ref{ss:5.4}} as follows: for $p\equiv1 \mod 4$
$$\binom{\frac{p-1}2}{\frac{p-1}4}\equiv (-4)^{\frac{p-1}4}(a+b\sqrt{-1}) \mod p^2, \quad \text{where}~ p=a^2+b^2 ~\text{with}~ a\equiv 1\mod 4.$$ {Here $\sqrt{-1}$ is a fixed unit in $\Z_p$ of order $4$, and $b$ is chosen so that $a+b\sqrt{-1}$ is a unit in $\Z_p$.}  Such a congruence, which is stronger than what can be predicted by the theory of formal group, is called a \emph{supercongruence}. Examples of supercongruences include \eqref{eq:18} and \eqref{eq:22} above.  Here, we focus on the cases with complex multiplication. Using properties of the $p$-adic Gamma function, Coster \cite{Coster} was able to extend the above result to \begin{equation}\label{eq:Coster}\binom{\frac{p^r-1}2}{\frac{p^r-1}4}/\binom{\frac{p^{r-1}-1}2}{\frac{p^{r-1}-1}4}\equiv (-4)^{\frac{p^{r-1}(p-1)}4}(a+b\sqrt{-1}) \mod p^{2r}, \quad \text{with } a, b \text{ as above}.\end{equation} Coster and van Hamme have the following further generalization. We  use $P_k(x)$ to denote ${}_2F_1\left[{-k, 1+k \atop 1}; \frac{1-x}2 \right]$,
which is a degree-$k$ polynomial known as the $k$th Legendre polynomial.

 \begin{theorem}[Coster and van Hamme, \cite{CosterVanHamme}]\label{CosterVanHammeTheorem} Let $d$ be a square-free positive integer, and $K = \Q(\sqrt {-d})$. Suppose that the elliptic curve
\[\mathcal{E}: Y^2=X(X^2+AX+B)\]   over $K$ has complex multiplication by an order of the ring of integers of $K$. Let $p$ be a prime split in $K$ such that $A, B$, and the square roots of $\Delta=A^2-4B$ can be embedded in $\Z_p^\times$. Then there exists $\a_{\mathcal E,p} \in \Z_p$ such that  $$\displaystyle \lim_{r \rightarrow \infty}P_{\frac{p^r-1}{2}}\left(\frac{A}{\sqrt{\Delta}}\right) /P_{\frac{p^{r-1}-1}{2}}\left(\frac{A}{\sqrt{\Delta}}\right)  = \a_{\mathcal E,p}.$$  Moreover, \begin{equation}
P_{\frac{mp^r-1}{2}}\left(\frac{A}{\sqrt{\Delta}}\right) \equiv
\a_{\mathcal E,p} \cdot
P_{\frac{mp^{r-1}-1}{2}}\left(\frac{A}{\sqrt{\Delta}}\right) \mod p^{2r}~~~\forall  \text{ odd } m\ge 1.
\end{equation}
\end{theorem}

In the heart of their proof lies  a  special Frobenius lifting that commutes with the endomorphism ring $\frak R$ of $\mathcal E$. This Frobenius arises from a $p$-isogeny sending  $t=-\frac{X}Y$ to  a degree-$p$ rational function of $t$. Using the above theorem and Clausen's formula \eqref{eq:clausen}, one has
\begin{theorem}[Kibelbek, Long, Moss, Sheller, Yuan, \cite{KLMSY}]\label{thm:KLMSY}
Let $\lambda\neq 1$ be an algebraic number such that $\tilde E_\lambda: y^2=(x-1)(x^2-\frac{1}{1-\l})$ has
complex multiplication. Let $p$ be a prime {such that $1-\l$ can be embedded in $\Z_p^\times$}  and $\tilde E_\l$  has   good
reduction modulo $p{\mathbb Z}_p$.
 Then
 $${}_3F_2\left[{\frac12, \frac12, \frac12 \atop 1, 1}; \l \right]_{\frac{p-1}2} =
\sum_{k=0}^{\frac{p-1}{2}}
\left(\frac{(\frac 12)_k}{k!}\right)^3\l^k\equiv \left( \frac{ 1-\l}{p} \right)\alpha_{p,\lambda}{}^2 \mod {p^2}$$
where $ \left( \frac{1-\l}{p} \right )$ is the Legendre symbol as before, $\alpha_{p,\lambda}$ is the unit root of
$X^2-[p+1-\#(\tilde E_\l/{\mathbb F}_p)]X+p$ if  $\tilde E_\l$ is ordinary at $p$; and
$\alpha_{p,\lambda}=0$ if $\tilde E_\lambda$ is supersingular at $p$.
\end{theorem}
As a corollary, one can deduce the following type of results that are difficult to prove directly.
\begin{cor}
For  all primes $p>3$, we have
 \begin{equation}\label{cor:4}
   \sum_{i=1}^{\frac{p-1}2} \binom{2i}{i}^3 \sum_{j=1}^{i} \frac{1}{i+j} \equiv 0 \mod{p}.
 \end{equation}
\end{cor}

Moreover, numerical data suggest that for any  ordinary prime  $p>3$ satisfying the above condition and any $m,n\ge 1$ and $m$ odd
\begin{equation}
{}_3F_2\left[{\frac12, \frac12, \frac12 \atop 1, 1}; \l \right]_{\frac{mp^n-1}2}\equiv
\left( \frac{{1-\l}}{p} \right)\alpha_{p,\lambda}{}^2\,
{}_3F_2\left[{\frac12, \frac12, \frac12 \atop 1, 1}; \l \right]_{\frac{mp^{n-1}-1}2} \mod
 {p^{3n}}.\end{equation}
In \cite{Villegas}, Rodriguez-Villegas made several conjectures on supercongruences  relating truncated hypergeometric series to Hecke eigenforms. Many of his conjectures have been proved and often the proofs rely on results like \eqref{cor:4}.

\subsection{ASD congruences and Fuchsian ordinary differential equations}\label{ss:6}
Ap\'ery numbers play an important role in transcendence. Using them Ap\'ery showed the irrationality of zeta values $\zeta(2)$ and $\zeta(3)$. Here, we will focus on the sequence $$A(n)=\sum_{k=0}^n \binom{n}{k}^2\binom{n+k}{k}=\,{}_3F_2\left[{-n, -n, n+1 \atop 1, 1}; 1 \right]$$
that  is related to $\zeta(2)$. This sequence satisfies several surprising ASD-type congruences. In \cite{Beu2}, Beukers showed that for a prime $p>3$,
\begin{equation}\label{eq:Beukers}
  A(mp^n-1)\equiv A(mp^{n-1}-1) \mod p^{3n}, \quad \forall m,n\ge 1.
\end{equation}
In \cite{SB}, Stienstra and Beukers proved that for  prime $p>3$, {$n \ge 1$ and odd $m \ge 1$, }
\begin{equation}\label{eq:SB1}
A(\frac{mp^n-1}2)-a_pA(\frac{mp^{n-1}-1}2)+\left (\frac{-1}p\right)p^2 A(\frac{mp^{n-2}-1}2)\equiv 0 \mod p^n,
\end{equation}  where $a_p$ is the $p$th Fourier coefficient of  $\eta(4z)^6$ appeared before. Stienstra-Beukers conjectured that \eqref{eq:SB1}  holds mod $p^{2n}$. When $m=n=1$, the conjecture was proved by Ishikawa \cite{Ishikawa} and Ahlgren \cite{Ahlgren}.

Now we briefly explain the geometry behind the above results using the viewpoint of \cite{Beukers}. It is well-known that $\sum_{n\ge 0} A(n)t^n$ satisfies a second order ordinary differential equation (ODE)
\begin{equation}\label{eq:G1(5)}
t(t^2-11t-1)\frac{d^2F(t)}{dt^2}+(3t^2-22t-1)\frac{dF(t)}{dt}+(t-3)F(t)=0.
\end{equation}In \cite{Hon71}, Honda asked  ``what algebraic differential equations `yield' formal groups that
are integral for almost all primes". Based on \cite{Dwork-p-cycle} Honda in \cite{Hon72} constructed formal power series satisfying linear ODEs that yield formal groups related to the Fermat curves. The coefficients of these power series are similar to the sequence used in \eqref{eq:Coster}.
Honda's question concerns when a given linear Fuchsian (i.e. with regular singularities only) ODE $Lf=0$ arises from a Picard-Fuchs equation.  A deep theorem of Katz says that Picard-Fuchs equations are globally nilpotent, see \cite{Katz70} for terminology and details. Let $\Sigma$ denote the set of regular singularities of $Lf=0$. In particular, we restrict ourselves to second order ODE.
Recall that given such a Fuchsian ODE, one can construct a \emph{monodromy representation}  $$\rho: \pi_1(\C P^1\setminus \Sigma, u_0) \rightarrow GL_2(\C)$$ of the fundamental group $\pi_1(\C P^1\setminus \Sigma, u_0)$, where $u_0$ is any nonsingular point on the base curve \cite{Yoshida}.  In particular, if $\text{Im} \rho$ can be embedded to $SL_2(\mathbb R)$,
then the image $\text{Im} \rho$ is a Fuchsian group \cite{shim1}. In this case,  the local holomorphic solution of $Lf=0$ (under suitable assumption) is a weight-1 automorphic form of the Fuchsian group $\text{Im} \rho$, like \eqref{eq:3}, see \cite{Stiller}. Upon knowing the local monodromy matrix at each singular point, one can tell whether $\text{Im} \rho$ has cusps or not. One can further ask whether $\text{Im} \rho$ is \emph{commensurable} with $SL_2(\Z)$. When the ODE has only 3 singularities, i.e. being a hypergeometric differential equation, the answer is known due to classification of arithmetic triangular groups \cite{Tak}. Next in line are second order ODEs with 4 singularities, like \eqref{eq:G1(5)}. For \eqref{eq:G1(5)}, the monodromy group $\text{Im} \rho$ is  isomorphic to $\Gamma^1(5)$ mentioned in \S2. The modular curve for $\G^1(5)$ has genus 0 and we can pick as a Hauptmodul $t=\frac{E_2}{E_1}=q^{1/5}+\cdots$, where $E_1,E_2$ are weight-3 Eisenstein series for $\G^1(5)$ that appeared in \S \ref{ss:mf}. Like \eqref{eq:3}, when we specify $t$ as a modular function for $\G^1(5)$, $\sum_{n\ge 0} A(n) t^n$ is a weight-1 modular form and  $\frac{dt}{dq^{1/5}}$ is a weight-2 modular form. Putting together
one can verify that
$$\sum_{n\ge 0} A(n) t^n\frac{dt}t=E_1\frac{dq_5}{q_5},\quad q_5=q^{1/5}$$ where $E_1(z)$ is a Hecke eigenform for all Hecke operators $T_p$ when $p>5$ with 1 as an eigenvalue.
Thus Beukers' result \eqref{eq:Beukers} modulo $p^n$ can be explained via Proposition \ref{prop:changevariable-CFGL} { and Theorem \ref{thm:iso-Hilbert}} as we are expressing  $E_1\frac{dq_5}{q_5}$ in terms of $t$, which is another local uniformizer at infinity.   Meanwhile, we emphasize that the change of variable as described in Proposition \ref{prop:changevariable-CFGL} does not preserve supercongruences in general.
See \cite{CCS} for some conjecture similar to \eqref{eq:Beukers} satisfied by Ap\'ery-like sequences.

To see \eqref{eq:SB1} we adopt a similar viewpoint.  Let $t_2=\sqrt{t}$. Thus $t_2E_1=\sqrt{E_1E_2}$ is the weight-3 cusp form mentioned in \S2. Then
$$ t_2E_1\frac{dq_{10}}{q_{10}}=\sum_{n\ge 0} A(n)t_2^{2n+1}\frac{dt_2}{t_2}, \quad q_{10}=q^{1/10}.$$

Motivated by results of   Beukers on Ap\'ery numbers, Zagier \cite{Zagier} did an extensive computer search for  rational numbers $a,b,\lambda$ such that the  differential equation
\begin{equation}\label{eq:1}
  (t(t^2+at+b)F'(t))'+(t-\l)F(t)=0
\end{equation}
has a solution in $\Z[[t]]$. Among his findings are cases where the monodromy groups are finite and ODEs  are equivalent to those of {hypergeometric series.} The most interesting sporadic cases have 4 genuine singularities and infinite monodromy groups, in which the monodromy group is isomorphic to an index-12 subgroup of $\SL$.  {Zagier conjectured that there are no other $(a,b,\l)$ giving rise to ODE with  solutions in $\Z[[t]]$. For the sporadic cases, his} conjecture is actually implied by an earlier conjecture by Chudnovsky-Chudnovsky on special globally nilpotent  Lam\'e  ODEs with 4 singularities, see \cite[Conjecture 2]{CC}.

\section{Weight-$k$ ASD congruences for noncongruence modular forms}\label{ss:7}
As alluded to earlier, congruences that are stronger than what commutative formal group laws  predict are harder to achieve, and their existence is usually due to extra symmetry like complex multiplication and/or special choice of the local uniformizer as well as the Frobenius lifting. However, in the realm of noncongruence modular forms of weight $k > 2$, one can always achieve supercongruence. This is due to the very special local uniformizer $q=e^{2\pi i z}$, the Frobenius lifting $q\mapsto q^p$ and their relation with the Gauss-Manin connection. For details, see
Katz \cite{Kat}.

\subsection{ASD congruences for weakly holomorphic modular forms} Let $\G$ be a finite index subgroup of $SL_2(\mathbb Z)$. Assume that the modular curve $X_\G$ has a model over $\Q$, the cusp at $\infty$ is a $\Q$-rational point with cusp width $m$. The completion of the local ring $\O_{X_\G, \infty}$ at $\infty$ is isomorphic to $\Q[[t]]$ for some $t$ satisfying $\delta t^m = q$ with $\delta \in \Q^\times$. Recall that there is an integer $M$, divisible by the widths of the cusps of $\G$ and the primes $p$ where the cusps are no longer distinct under reduction of $X_\G$ at $p$, such that for any integer $k \ge 1$, the space
of weight-$k$ modular forms for $\G$ has a basis whose $t$-expansions have coefficients in $\Z[\frac1M]$.
Given a subring $R$ of $\C$ such that $6M$ is invertible in $R$, {the modular curve  $X_{\G}$ has a model defined over $R$ and $X_\G(R)$ contains the cusp infinity,}  let $M_k^{wk}(\G,R)$ be the set of weakly holomorphic (i.e., holomorphic on the upper half-plane and meromorphic at {all the cusps}) weight $k$ modular forms for $\G$ whose $t$-expansions at $\infty$ have coefficients in $R$. Denote by $S_k^{wk}(\G,R)$ the submodule of  functions in $M_k^{wk}(\G,R)$ with vanishing constant {terms at all the cusps.
{The Fourier coefficients $a_n(f, c)$ of $f \in M_k^{wk}(\G,R)$ at any cusp $c$ of $\G$ are integral over $R$.} A modular form $f\in M_{k}^{wk}(\G,R)$ is called \emph{weakly exact} if {$n^{-(k-1)}a_n(f,c)$ is integral over $R$ for each $n<0$ and each cusp $c$ of $\G$.}
Write $S_k^{wk-ex}(\G,R)$ (resp. $M_k^{wk-ex}(\G,R)$) for the collection of weakly exact forms in $S_k^{wk}(\G,R)$ (resp. $M_k^{wk}(\G,R)$). Assume $k \ge 2$.  The linear map $D=q\frac{d}{dq}$ iterated $k-1$ times maps $M_{2-k}^{wk}(\G,\C)$ into $S_k^{wk}(\G,\C)$ {by Bol's identity \cite{Bol}}.
{Furthermore,  given $f \in M_{2-k}^{wk}(\G, R)$, at each cusp $c$ of $\G$ with cusp width $m_c$,  we have $a_n(D^{k-1} f, c) = (\frac{1}{m_c})^{k-1} n^{k-1} a_n(f, c)$, which shows that $D^{k-1} f $ lies in
$S_k^{wk-ex}(\G, R)$ since $1/m_c \in R$ by assumption.}

Using geometric interpretations of weakly exact forms, Kazalicki and Scholl in \cite{KS13} identified the quotient space $$DR(\G,R,k):=\frac{S_k^{wk-ex}(\G,R)}{D^{k-1}(M_{2-k}^{wk}(\G,R))}$$with the $p$-adic de Rham space which plays a key role in the proof of the Scholl congruence (\ref{schollcongruence}) {\cite[\S 3]{Sch85b}}.
Similarly define
$$DR^*(\G,R,k):=\frac{M_k^{wk-ex}(\G,R)}{D^{k-1}(M_{2-k}^{wk}(\G,R))}.$$
As there are no holomorphic forms of negative weight, $S_k(\G, R)$ is contained in $DR(\G, R, k)$ and $M_k(\G, R)$ in $DR^*(\G, R, k)$. In fact, the following two short exact sequences {hold due to Serre  duality:
$$0\rightarrow S_k(\G,R)\rightarrow DR(\G,R,k)\rightarrow S_k(\G,R)^{\vee}\rightarrow 0$$ and
$$0\rightarrow M_k(\G,R)\rightarrow DR^*(\G,R,k)\rightarrow S_k(\G,R)^{\vee}\rightarrow 0,$$ where
  $S_k(\G,R)^{\vee}$ is the $R$-linear dual of $S_k(\G,R)$.}
 Suppose that $S_k(\G, R)$ has $R$-rank $d = d(k)$. The above exact sequences imply that $DR(\G,R,k)$ and $DR^*(\G,R,k)$ are locally free $R$-modules of rank $2d$ and $2d$ plus the rank $d'(k)$ of the weight $k$ Eisenstein series in $M_k(\G, R)$, respectively.

It was shown in  \cite{A-SD} and \cite{Sch85b}  that there is a positive integer $M$ such that the $d$-dimensional space $S_k(\G)$ has a basis consisting of functions whose $t$-expansions { at the cusp $\infty$} have coefficients in $R=\mathbb Z[\frac{1}{M}]$ and their $q$-expansions are integral over $R$. For any prime $p\nmid M$, we have $R$ contained in $\Z_p$  (embedded  in $\C$).  Further, for $p>k-2$, there is an endomorphism $\phi_p$ on $DR^*(\G, \Z_p, k)$ leaving $DR(\G, \Z_p, k)$ invariant, arising from  the  Frobenius lifting originated in the map $q\mapsto q^p$ on the Tate curves.   Its characteristic polynomial $H_p(T)$
 on $DR(\G, \Z_p, k)$ lies in  $\Z[T]$, and coincides with the characteristic polynomial of the geometric Frobenius at $p$ under the $2d$-dimensional Scholl representations recalled in \S \ref{ss:mf}.    Thus all roots of $H_p(T)$ are algebraic integers with the same absolute value $p^{(k-1)/2}$, and the non-real complex roots can be paired off as $\{\alpha, p^{k-1}/\alpha\}$.  The characteristic polynomial of $\phi_p$ on the quotient $DR^*(\G, \Z_p, k)/DR(\G, \Z_p, k)$ also lies in $\Z[T]$, {with all roots of absolute value ${ p^{k-1}}$} (cf. \cite[pp. 75]{Sch85b}).   Kazalicki and Scholl proved in \cite{KS13} that the congruences (\ref{schollcongruence}) satisfied by the cusp forms in $S_k(\G, \Z_p)$ also hold for weakly holomorphic forms in $M_k^{wk-ex}(\G, \Z_p)$, although with weaker moduli.  {Since to study the behavior at another cusp amounts to replacing $\G$ by a conjugate, we only consider the cusp at infinity.}

\begin{theorem}[\cite{KS13}]\label{thm:KS} Let $p>k-1$ be a prime {such that $p\nmid 6M$}. Suppose $f = \sum a_n(f) q^{n/m}$ in { $M_k^{wk-ex}(\G, \Z_p)$} is annihilated by $h(\phi_p)$ for some polynomial $h(T)=\sum_{j=0}^{r}A_jT^j\in \Z[T]$ dividing the characteristic polynomial of $\phi_p$ on $DR^*(\G, \Z_p,k)$. Then
{$$\sum_{j=0}^{r} p^{(k-1)j}A_ja_{n/p^j}(f)\equiv 0 \mod p^{(k-1)\ord_p n}, \qquad \text{ for~ all}~ n \ge 1.$$}
\end{theorem}

\begin{example}The space $S_{12}^{wk-ex}(\SL,\Z)$ is a $\Z$-module spanned by $\Delta(z) = \sum_{n \ge 1} \tau_nq^n$ and
$$g(z) := E_4(z)^6/\Delta(z)-1464E_4(z)^3=q^{-1}-1432236q+51123200q^2+39826861650q^3+\cdots.$$
As is well-known, the Fourier coefficients $\tau_n$ of $\Delta$ satisfy the recursion $$\tau_{np} - \tau_p\tau_n + p^{11}\tau_{n/p} = 0 \quad \text{for~all}~ n\ge 1 ~\text{and~all~primes}~ p.$$  Kazalicki and Scholl showed that,
 for every prime $p\ge 11$, the Fourier coefficients $a_n(g)$ of $g$ satisfy
$$a_{np}(g)-\tau_pa_n(g)+p^{11}a_{n/p}(g)\equiv 0 \mod p^{11\ord_p n}\quad \text{for~all}~ n\ge 1.$$

\end{example}

\subsection{ASD congruences for  cusp forms}
The stronger congruences (\ref{schollcongruence}) satisfied by cusp forms $f\in S_k(\G, \Z_p)$ established by Scholl in \cite{Sch85b} resulted from the fact that $\phi_p$ on $DR(\G, \Z_p, k)$ actually sends $S_k(\G, \Z_p)$ into $p^{k-1}DR(\G, \Z_p, k)$. The extra multiple $p^{k-1}$ accounts for the higher exponent
in the moduli. More precisely, the Fourier coefficients $a_n(f)$ of $f$ at the cusp $\infty$ satisfy the congruence
$$\sum_{j=0}^{2d} p^{(k-1)j}A_ja_{n/p^j}(f)\equiv 0 \mod p^{(k-1)(1+\ord_p n)}, \qquad \text{ for ~all}~ n \ge 1,$$where $d$ is the dimension of $S_k(\G)$ and $H_p(T) = \sum_{j=0}^{2d} A_jT^j$ is the characteristic polynomial of $\phi_p$.

 J. Kibelbek in his thesis \cite{Kibelbek}  gave an interpretation of the above $(2d+1)$-term ASD congruences in terms of $d$-CFGLs for the case $k=2$.

 In general, due to the lack of effective Hecke operators, one does not know how to decompose the { $2d$-dimensional} Scholl representations into a sum of 2-dimensional subrepresentations, as what happened for Deligne representations for congruence groups.  However, if extra symmetries are present, then sometimes they can be used to break the Scholl representations. Accordingly, one can decompose Scholl congruences into $3$-term ASD congruences, resembling the congruence case. We demonstrate below a few cases where extra symmetries are used to obtain 3-term ASD congruences for almost all primes and (semi-)fixed basis. More examples can be found in {\cite{LLY05, ALL05, Long08, L5, HLV, ALLL}.}

\subsection{Examples}
Let $\G_3$ be an index-3 subgroup of $\G^1(5)$ such that $t_3:=\sqrt[3]{t} = \sqrt[3]{E_2/E_1}$ is a Hauptmodul of the modular curve $X_{\G_3}$.
The space $S_3(\G_3)$ is spanned by $g_1=E_1t_3$ and $g_2=E_1t_3^2$,  both in $\Z[1/3][[q^{1/15}]]$.
\begin{theorem}[Li, Long, Yang \cite{LLY05}]For any prime $p>3$, the coefficients of $g_1\pm \sqrt{-1}g_2=\sum_{n \ge 1} a_{\pm}(n)q^{n/15}$ satisfy
$$a_{\pm}(np^r)-b_{\pm}(p)a_{\pm}(np^{r-1})+\chi_{-3}(p)p^2 a_{\pm}(np^{r-2}) \equiv 0 \mod p^{2r}, \quad \text{for}~ r,n\ge 1,$$ where $\chi_{-3}$  is the quadratic character associated to $\Q(\sqrt {-3})$, and $g_{\pm}(z):=\sum _{n\ge 1}b_{\pm}(n)q^n,$ $q=e^{2\pi i z}$ are two weight-$3$ normalized congruence cuspidal newforms of level $27$ and character $\chi_{-3}$ {whose Fourier expansions start with
\begin{eqnarray*}
g_+(z)&:=&q-3\sqrt{-1}q^2-5q^4+3\sqrt{-1}q^5+5q^7+3\sqrt{-1}q^8+9q^{10}+\cdots\\
g_-(z)&:=&q+3\sqrt{-1}q^2-5q^4-3\sqrt{-1}q^5+5q^7-3\sqrt{-1}q^8+9q^{10}+\cdots.
\end{eqnarray*}}
\end{theorem}
Next we consider the index-$4$ subgroup $\G_4$ of $\G^1(5)$ defined similarly using $t_4 = \sqrt[4]{E_2/E_1}$.
Let $S$ be the space generated by the two weight-$3$ cusp forms $h_1=E_1t_4$ and $h_3=E_1 t_4^3$ for $\G_4$.  There is a compatible family of $4$-dimensional sub-Scholl-representations of $G_\Q$ attached to $S$. Consider the following four weight-$3$ congruence cuspforms defined using the $\eta$ function: $$f_1(z)= \frac{\eta(2z)^{12}}{\eta(z)\eta(4z)^5}=\sum_{n\ge 1}a_1(n)q^{n/8}, \qquad \quad  f_3(z)=\eta(z)^5\eta(4z)=\sum_{n\ge 1}a_3(n)q^{n/8},$$ $$f_5(z)= \frac{\eta(2z)^{12}}{\eta(z)^5\eta(4z)}=\sum_{n\ge 1}a_5(n)q^{n/8},~~~~~~ \qquad ~\text{and}~\quad  f_7(z)=\eta(z)\eta(4z)^5=\sum_{n\ge
1}a_7(n)q^{n/8}.$$
\noindent Their linear combination \begin{equation} \label{eq:f'}f=f(z) = f_1(z)
+ 4f_5(z) + 2\sqrt{-2}(f_3(z) - 4f_7(z)) = \sum_{n \ge
1}a(n)q^{n/8}\end{equation} is an eigenform of
the Hecke operators at odd primes and $f(8z)$ has level 256, weight 3, and quadratic character
$\chi_{-4}$ associated to {$\Q(\sqrt{-1})$}.

\begin{theorem}[Atkin, Li, Long \cite{ALL05}]
For each  prime odd $p>2$, the space {$S =\langle h_1, h_3\rangle$} has a basis
depending on the residue of $p \mod 8$ satisfying the ASD
congruence \eqref{ASDgeneral} at $p$ as follows.
\begin{enumerate}
\item If $p \equiv 1 \mod 8$, then both $h_1$ and $h_3$ satisfy \eqref{ASDgeneral} with the same $A_p =\text{sgn}(p)a_1(p)$  and $\mu_p = 1$, where $\text{sgn}(p) = \pm 1 \equiv 2^{(p-1)/4} \mod p$;
\item If $p \equiv 5 \mod 8$, then $h_1$ (resp. $h_3$) satisfies \eqref{ASDgeneral} with $A_p= 4\sqrt{-1}a_5(p)$ (resp. $-4\sqrt{-1}a_5(p)$)  and $\mu_p = -1$;
 \item If $p \equiv 3 \mod 8$, then $h_1\pm h_3$  satisfy \eqref{ASDgeneral} with $A_p= \pm 2\sqrt{-2}a_3(p)$  respectively, and $\mu_p = -1$;
  \item If $p \equiv 7 \mod 8$, then $h_1\pm \sqrt {-1}h_3$  satisfy \eqref{ASDgeneral} with $A_p= \mp 8\sqrt{-2}a_7(p)$ respectively, and $\mu_p = -1$.
  \end{enumerate}
Here, for $j \in \{ 1, 3, 5, 7\}$, $a_j(p)$ is the $p$th Fourier coefficient of the congruence form $f_j$ defined above.
\end{theorem}

In \cite{Sch85b} Scholl showed that if half of the roots of the characteristic polynomial $H_p(T)$ of $\phi_p$ are distinct $p$-adic units, then one can find a basis for $S_k(\G)$ satisfying the $3$-term ASD congruence. In his thesis J. Kibelbek constructed examples to show that the $3$-term ASD congruence does not always hold  \cite{Kib12}.  A similar example was given in
\cite{KS13}. The example below is due to Kibelbek.

Consider the genus $2$ hyperelliptic curve $X$ over $\Q$ with an affine equation  $y^2 = x^5 + 2$. {The Jacobian of this this curve admits complex multiplication.} By Bely{\u\i}'s theorem, $X \simeq X_\G$ for a finite index subgroup $\G$ of $SL_2(\Z)$. {One isomorphism is given by letting $x=-\sqrt[5]{2\l(z)}, \, y=\sqrt{2-2\l(z)}$} where $\l(z)$ is the modular lambda function given in \S\ref{ss:4.1}.  In this realization, the cusp at infinity has cusp width $10$.  On $X$ there are two linearly independent holomorphic differentials $\frac{dx}{2y}$ and $x\frac{dx}{2y}$.  Using the $q$-expansion of $\l(z)$, they can be can rewritten, up to a normalization by a constant multiple, as $f_1 \frac{dq^{1/10}}{q^{1/10}}$ and $ f_2 \frac{dq^{1/10}}{q^{1/10}},$ respectively, where

\begin{eqnarray*}
f_1 &=& q^{1/10} - \frac{8}{5}q^{6/10} - \frac{108}{5^2}q^{11/10} + \frac{768}{5^3} q^{16/10} + \frac{3374}{5^4}q^{21/10} + \cdots= \sum_{n \ge 1} a_n(f_1)q^{n/10}
\end{eqnarray*} and
\begin{eqnarray*}
f_2 &=&q^{2/10} - \frac{16}{5}q^{7/10} + \frac{48}{5^2}q^{12/10} + \frac{64}{5^3} q^{17/10} + \frac{724}{5^4}q^{22/10} + \cdots= \sum_{n \ge 1} a_n(f_2)q^{n/10}.
\end{eqnarray*} They generate the space $S_2(\G)$ and the module $S_2(\G, \mathbb Z[\frac{1}{5}])$.

The {$\ell$-adic} representations attached to $S_2(\G)$ are the dual of the Tate module of the Jacobian of $X_{\Gamma}$. At primes $p \equiv 1,4 \mod 5$, {$f_1$ and $f_2$ each satisfy a  $3$-term ASD congruence at $p$, as a consequence of the CM structure. For primes $p \equiv 2, 3 \mod 5$,   $H_p(T) = T^4 + p^2$.}  At such primes both $f_1$ and $f_2$ }and consequently their nontrivial linear combinations $g$ satisfy the Scholl congruences $$a_{mp^{n+2}}(g) + p^2a_{mp^{n-2}}(g)\equiv 0 \mod p^{n+1}, \quad \quad \forall n,m\ge 1, $$ but no $3$-term congruences relating $a_{mp^n}(g), a_{mp^{n+1}}(g)$ and $a_{mp^{n+2}}(g)$ would hold since $p \equiv 2, 3 \mod 5$ and $a_n(f_i)\ne 0$ only when $n\equiv i \mod 5$ for $i = 1, 2$.
 On the other hand, as explained below, the module of weakly holomorphic forms $DR(\Gamma, \Z_p[\sqrt {2p}], 2)$ does contain a basis satisfying 3-term ASD-type congruences as Theorem 13 for primes $p \ne 5$ since $H_p(T)$ factors into a product or two degree-$2$ polynomials over $\Z_p[\sqrt {2p}]$.

The differentials of the second kind $x^2\frac{dx}{2y}$ and $x^3\frac{dx}{2y}$ on $X$ give rise to two weakly holomorphic cusp forms
\begin{eqnarray*}
f_3 &=& q^{3/10} - \frac{24}{5} q^{8/10} + \frac{268}{5^2}q^{13/10} - \frac{2624}{5^3}q^{18/10} + \frac{24714}{5^4}q^{23/10}  + \cdots= \sum_{n \ge 1} a_n(f_3)q^{n/10},
\end{eqnarray*} and
\begin{eqnarray*}
f_4 &=&  q^{4/10} - \frac{32}{5}q^{9/10} + \frac{552}{5^2}q^{14/10} - \frac{7808}{5^3}q^{19/10} + \frac{97104}{5^4}q^{24/10} + \cdots=\sum_{n \ge 1} a_n(f_4)q^{n/10}.
\end{eqnarray*}
Note that $f_3$ and $f_4$ are holomorphic at $\infty$ but have poles at other cusps. The four forms $f_1, f_2, f_3, f_4$ together span the space $DR(\G, \Z[\frac{1}{5}], 2)$. It is straightforward to check that ASD congruences at primes $p \ne 5$ with weaker moduli as in Theorem \ref{thm:KS} are satisfied by four linearly independent forms in  {$DR(\G, \Z_p[\sqrt {2p}], 2)$.}

\section{An application of ASD congruences}\label{ss:app}
Congruence forms with algebraic Fourier coefficients are known to have bounded denominators. A folklore conjecture asserts that a cusp form for a finite index subgroup of $SL_2(\Z)$ with algebraic Fourier coefficients is a congruence form if and only if its Fourier coefficients have bounded denominators. When the space $S_k(\G)$ is $1$-dimensional, as discussed before, the ASD congruences hold and the associated Scholl's representations are modular. Using these facts, we established in \cite{LL12} the conjecture for the $1$-dimensional case.

\begin{theorem}[Li and Long, \cite{LL12}]\label{thm:LL} { Suppose that the modular curve $X_\G$ of $\G$ has a model defined over $\Q$ so that the cusp at $\infty$ is $\Q$-rational { with cusp width $m$},  $k\ge 2$ and $S_k(\G)$ is 1-dimensional.} Then a form $f=\sum a_n(f)q^{n/m}$ in $S_k(\G)$ with Fourier coefficients in $\Q$ has bounded denominators if and only if it is a congruence cusp form.
\end{theorem}
We outline the idea of the proof. As $f$ has bounded denominators, without loss of generality we may assume that $a_n(f)\in \Z$. Selberg \cite{selberg65} proved that the Fourier coefficients of $f$ satisfy the following bound for some constant $C$:
\begin{equation}\label{eq:bn}
  |a_n(f)|<Cn^{k/2-1/5}, \quad \forall n\ge 1.
\end{equation}Meanwhile,
{the compatible family of $\ell$-adic Scholl representations attached to $S_k(\G)$ is modular}, that is, the associated $L$-function coincides with the $L$-function of a cuspidal congruence Hecke eigenform $g=\sum_{n \ge 1} b_nq^n$  {of weight $k$, level $N$, character $\chi$, and integer coefficients. The ASD congruence established by Scholl says that there is an integer $M$ such that for $p\nmid M$
$$a_{np^s}(f)-b_pa_{np^{s-1}}(f)+\chi(p)p^{k-1}a_{np^{s-2}}(f)\equiv 0 \mod p^{(k-1)s}, \quad \forall n\ge 1.$$ The left-hand side is an integer bounded by a fixed constant multiple of $(np^s)^{k/2-1/5}$ for all $n$ and $s$. Thus for fixed $n$ and $s$ large enough, the congruence becomes an equality. This implies that, after twisting by a suitable multiplicative character $\psi$ of $\Z$ to get rid of multiples of small primes, $f_\psi(z)=\sum_{n \ge 1} \psi(n)a_n(f)q^n$ and  $g_\psi(z)=\sum_{\ge 1} \psi(n)b_nq^n$ agree {up to a nonzero scalar multiple}. Since twisting by characters preserves the congruence/noncongruence property of modular forms, we conclude that $f$, which has bounded denominators, has to be a congruence form.

\section{Another type of congruences}Let $f$ be a meromorphic modular function for {$\SL$} with Fourier coefficients $c_n(f)\in \Z$,  i.e.,
$\displaystyle f=\sum_{n> -\infty} c_n(f)q^n.$ Fix a prime $p$. For a positive integer $m$ {such that $c_{p^m}(f) \ne 0$}, let $t_m(f,n)=c_{np^m}(f)/c_{p^m}(f)$. Atkin had done extensive investigation  on the coefficients of  the modular $j$-function $j(z)=q^{-1}+744+196884q+\cdots$. In \cite{Atkin67}, he proved that for $p=11$,
\begin{equation}\label{eq:Atkin-11}c_{n11^m}(j)\equiv 0 \mod 11^m, \quad \forall m,n\ge 1,
\end{equation} extending results of Lehner for $p=5,7$.
\begin{conj}[Atkin and O'Brien  \cite{AOB} and Atkin \cite{Atkin68}]For any prime $p\neq 13$ {such that $c_{p^m}(j)\neq 0$ for all $m\ge 1$,}
\begin{equation}
t_m(j,np)-t_m(j,n)t_m(j,p)+p^{-1}t_m(j,n/p)\equiv 0 \mod 13^{m}, \quad \forall m, n\ge 1
\end{equation} and
\begin{equation}\label{eq:32}
  t_m(j,13n)-t_m(j,n)t_m(j,13)\equiv 0 \mod 13^{m}, \quad \forall m, n\ge 1.
\end{equation}
\end{conj}

For a given prime $p$,  the Atkin $U_p$-operator sends  $f=\sum a_nq^n$ to $U_p(f)=\sum a_{pn}q^n$. Thus,
$$\frac{U_p^m(j-744)}{c_{p^m}(j-744)}=\sum _{n\ge 1}t_m(j-744,n)q^n.$$ The above conjecture was proved by Koike \cite{Koike} and Katz \cite[\S 3.13]{Kat} by knowing that a repeated application of $U_{13}$ to $j-744$ leads to a single $13$-adic Hecke eigenform.  The other akin conjectures of Atkin for primes $p\le 23$ are similarly  established by Guerzhoy \cite{Gue1,Gue2}. Atkin's observations and $U_p$ operator play important roles in the development of $p$-adic modular forms \cite[et al.]{Serre, Dwork73, Kat, Hida86, Coleman}.

Suppose that $f=\sum_{n > -\infty} c_n(f)q^n \in \Z_p((q))$ and $c_1(f)$ is a $p$-adic unit. If $f$ satisfies a 2-term ASD congruence like \eqref{eq:2-term} at $p$, that is, there exists a $p$-adic unit $\a_p$ such that
$$c_{np^{m+1}} (f) \equiv \alpha_p c_{np^m}(f) \mod p^{m+1}  \quad \text{for~all}~ m \ge 0 ~~\text{and}~n\ge 1,$$
then $c_{p^m}(f)$ is a $p$-adic unit for all $m \ge 0$ so that $t_m(f, n)$ are in $\Z_p$. In the congruence above we can replace $\a_p$ by $c_{p^{m+1}}(f)/c_{p^m}(f)$, and the resulting congruence is nothing but the congruence of the form \eqref{eq:32}:
$$ t_m(f, pn) - t_m(f, n)t_m(f,p) \equiv 0 \quad \mod p^{m+1}.$$

 Kazalicki \cite{Kazalicki} observed that a congruence of type \eqref{eq:Atkin-11} for $p=2$ is satisfied by a family of noncongruence modular functions. Recently, the second author and Alyson Denies computed the following during {\sf Sage Day 46} for a weight-1 noncongruence form.  By Sebbar \cite{Sebbar},  the space of integral weight holomorphic modular forms for $\G_1(5)$  is a graded algebra {generated by two normalized weight-$1$ forms $f_1,f_2$ (\cite[pp. 302]{Sebbar}), whose zeros are located at the cusps only.} Let $f = \sqrt{f_1 f_2}$, which is
a weight-$1$ noncongruence  modular form. Our data suggest the following pattern:
\begin{conj}
  Let $p = 5$. Then for $f = \sqrt{f_1 f_2}$, $m\ge 1$ and odd  $n\ge 1$, we have {$c_{5^m}(f)\neq 0$ and}
  $$t_m(f,5n)\equiv t_m(f,n) \mod 5^{2m+4}.$$
\end{conj}

\def\cprime{$'$}
\providecommand{\bysame}{\leavevmode\hbox to3em{\hrulefill}\thinspace}
\providecommand{\MR}{\relax\ifhmode\unskip\space\fi MR }
\providecommand{\MRhref}[2]{%
  \href{http://www.ams.org/mathscinet-getitem?mr=#1}{#2}
}
\providecommand{\href}[2]{#2}

\end{document}